\documentclass[12pt]{article}

\usepackage{times}
\usepackage{tabu}
\usepackage{graphicx}
\usepackage{multirow}
\usepackage{url}
\usepackage{amssymb}
\usepackage{amsmath}
\usepackage{pifont} 
\usepackage{lineno}
\usepackage{xcolor}
\usepackage{tikz}
\usepackage{cancel}

\usepackage{algorithm}
\usepackage[noend]{algpseudocode}
\usepackage{booktabs}  
\usepackage{makecell} 

\newenvironment{sciabstract}{%
\begin{quote} \bf}
{\end{quote}}

\title{High-order expansion of Neural Ordinary Differential Equations flows}

\author
{Dario Izzo$^{1\dagger*}$ and Sebastien Origer$^{1*}$ and Giacomo Acciarini$^{1*}$\\ and Francesco Biscani$^{2}$\\
\\
\normalsize{$^{1}$Advanced Concepts Team, European Space Agency}\\
\normalsize{European Space Research and Technology Centre (ESTEC)} \\
\normalsize{Keplerlaan 1, 2201 AZ Noordwijk, The Netherlands}\\
\\
\normalsize{$^{2}$Skylon Dynamics}\\
\normalsize{Kurpfalzstrasse 4a, 75053 Gondelsheim, Germany} \\
\\
\normalsize{$^\dagger$To whom correspondence should be addressed; E-mail:  dario.izzo@esa.int} \\
\\
\normalsize{$^*$ These authors contributed equally to this work.} \\
}

\date{}

\begin{document} 


\baselineskip24pt


\maketitle 

\hyphenation{wide-spread}


\begin{sciabstract}
Neural networks are redefining our use of ordinary differential equations, both by enhancing their ability to model complex phenomena that challenge classical approaches and by enabling their interpretation as infinite-depth neural models.
However, their practical applicability remains limited by the lack of explainability of the resulting dynamics, which is governed by an opaque neural network. Beyond numerical simulations, existing analytical approaches are largely restricted to leverage first-order gradient information due to computational constraints.
We thus introduce a new approach utilizing high-order differentials through complex mathematical constructs we call Event Transition Tensors, and show its use in analyzing the neural dynamics up to a generic event manifold. 
The versatility of our approach is demonstrated by analyzing neural state-feedback systems, characterizing uncertainties in a data-driven prey-predator control model, and mapping landing trajectories in a three-body neural Hamiltonian ordinary differential equation. In all cases, the method enhances the rigor of the resulting analysis by describing the neural dynamics through explicit mathematical constructs. By leveraging and expanding the differentiable programming toolkit, these results contribute to a deeper understanding of event-triggered neural ordinary differential equations as models of complex systems and provide a valuable tool to  explain their behavior.

\end{sciabstract}

\section*{Introduction}
Ordinary differential equations (ODEs) provide a powerful framework for modeling a wide range of phenomena, from epidemiology and ecology to cosmology, engineering, chemistry, neuroscience, and meteorology. 
However, in each case, unmodeled and unknown effects can impact the analysis of the resulting dynamics, limiting the accuracy and scope of model-derived conclusions. 
When modeled or observed data of the system are available, they can be leveraged to partially or fully identify the system, leading to data-driven approaches for discovering or refining ODE models.
Over the past two decades, data-driven approaches to ODE discovery and refinement have evolved from symbolic regression techniques \cite{tsoulos2006solving, schmidt2009distilling, izzo2017differentiable} to include methodologies based on artificial neural models of the system dynamics \cite{lejarza2022data, greydanus2019hamiltonian, chen2018neural, izzo2024optimality}. Recently, Buzhardt et al. (2025) \cite{buzhardt2024relationship} have shown the equivalence of these methods to state space methods and Koopman operator-based methods, strengthening neural learning of ODEs as a generic and powerful approach. Notably the universal approximation property \cite{hornik1989multilayer} of neural networks, revisited by Calin \cite{calin2020deep}, suggests that any ODE can, in principle, be represented by a neural model. 
Interestingly, the reverse has also been explored: Chen et al. \cite{chen2018neural} showed how deep learning models such as residual networks, recurrent neural network decoders, and normalizing flows can themselves be interpreted as differential equations. 
Their work introduced the term \emph{NeuralODE}, which has since gained popularity.
In subsequent studies, including this manuscript, NeuralODE broadly refers to ordinary differential equations in which some or all of the dynamics are represented by an artificial neural network. 
Such systems are now appearing with increasing frequency in very diverse scientific contexts such as Hamiltonian dynamics \cite{greydanus2019hamiltonian}, the analysis of ecological and evolutionary time-series data \cite{bonnaffe2021neural}, optimal control problems~\cite{izzo2024optimality, bottcher2022ai}, stiff chemical kinetics \cite{kumar2023physics} as well as autonomous driving~\cite{giang2024conditional} and spintronic modelling \cite{chen2022forecasting}.

The integration of ODEs into machine learning pipelines has spurred the development of tools for computing NeuralODE sensitivities, tailored specifically for the machine learning community. 
Consequently, established techniques, such as the adjoint method \cite{griewank1994computational} and variational equations \cite{beutler2005variational, skokos2010numerical}, are being incorporated into popular ML frameworks to facilitate the training and analysis of NeuralODEs. 
In line with Griewank’s astute observation that “gradients are really cheap” \cite{griewank1994computational}, prior studies have primarily focused on leveraging gradient (i.e. first order) information of the flow. Griewank, as well as Berz, are also widely credited for laying the foundations for the efficient computation of high-order derivatives~\cite{berz1988differential, corliss1997high, griewank2008evaluating}, which have been only recently implemented also in a machine learning context~\cite{jax2018github}. 
However, despite these existing foundational contributions, the practical application of the high-order derivatives of the system flow has remained limited in a NeuralODE setting, particularly for multivariate cases where the computational costs involved grow significantly. We thus present an efficient method to expand the flow of NeuralODEs at higher orders, specifically on predefined differential event manifolds, and we use this information to work towards a certifiable behavior of the underlying neural system avoiding to resort to Monte Carlo simulations or other approximate numerical methods.
The proposed approach leverages efficient Taylor integration~\cite{biscani2021revisiting} of the high-order variational equations associated with the NeuralODE, providing deeper insight in the propagation of uncertainties and adversarial perturbations within these systems. Our findings support the positioning of NeuralODEs as reliable tools for exploring a wide range of phenomena that are not suitable for modelling and prediction with traditional methods.

\section*{Results}
\subsection*{Event Transition Tensors}
Consider a generic NeuralODE \cite{chen2018neural} of the form:

\begin{equation}
\label{eq:dyn}
\left\{
\begin{array}{l}
\dot{\pmb x} = \pmb f(\pmb x, \pmb p, \mathcal N_{\tilde{\pmb \theta}}(\pmb x), t) \\
\pmb x(0) = \pmb x_0
\end{array}
\text{,}
\right.
\end{equation}
where $\pmb x$ is the state, $\pmb f$ the dynamics and the highly parametrized right-hand side contains both parameters $\pmb p$ related to the physical quantities and ones that pertain to the neural model $\tilde{\pmb \theta}$. 
Here, the ensemble of all parameters is indicated as $\pmb \theta = [\pmb p, \tilde{\pmb \theta}]$ and their dimensionality as $|\pmb \theta|$. 
The solution to the above problem is a \emph{trajectory} $ \pmb x(t)$, which is guaranteed to exist unique under relatively mild conditions on $\pmb f$. 
A more interesting object, for the purpose of our reasoning, is the \emph{flow} of the NeuralODE, a function represented via the symbol $\pmb \Phi(t; \pmb x_0, \pmb \theta)$. 
The flow represents all possible trajectories generated by varying initial conditions and parameters. 
It is a rather complex function, at the core of our methods, which we represent through its high-order Taylor expansion (here the multi-index notation is used to manipulate multivariate Taylor polynomials and the index $j$ to indicate a specific component of the flow):
\begin{equation}
    \Phi_j(t; \pmb x_0, \pmb \theta) \approx \mathcal{P}_{\Phi_j}^k(\delta \pmb x_0, \delta \pmb \theta)= \overline x_{jf} + \sum_{|\alpha|=1}^k \dfrac{1}{\alpha!}\partial^\alpha \Phi_j(t)\bigg|_{({\pmb{\overline x}}_0, {\pmb{\overline \theta}})} [\delta \pmb x_0, \delta \pmb \theta]^\alpha\text{.}
\end{equation}
The symbol $\mathcal P^k$ is used to indicate a multivariate polynomial of order $k$. Following Park \& Scheeres (2006) \cite{park2006nonlinear}, all coefficients of the above expansion (at the various orders) are called \emph{State Transition Tensors} (STTs).
These bear a great significance in the analysis of dynamical systems in aerospace and celestial mechanics~\cite{bani2019exact, boone2021orbital}, but have never been considered in the context of NeuralODEs (with the exception of \cite{izzo2020stability} where the first hint on their value to study neurocontrolled dynamics was given). The STTs encapsulate all information about the high-order derivatives of the flow at a given time $t$, and their efficient computation is enabled by implementations of automated differentiation techniques dealing specifically with high-order flow expansions \cite{berz2002cosy, izzo2020dcgp, jax2018github}.
Furthermore, we introduce an event manifold $e(\pmb x, \pmb \theta, t)$ (here, an event is a condition that is triggered along a trajectory $\pmb x(t)$, such as the arrival at a target surface, or the acquisition of a minimum velocity).
Its corresponding trigger time $t^*$ is the solution to the equation:
$$
e(\pmb \Phi(t; \pmb x_0, \pmb \theta), \pmb \theta, t) = 0 \rightarrow t^*(\pmb x_0, \pmb \theta)
$$
We may then introduce the \emph{event map} $\pmb F$ defined as the flow of the dynamical system at the event trigger time:
\begin{equation}
\label{eq:event_map}
\pmb F(\pmb x_0, \pmb p)=\pmb \Phi(t^*; \pmb x_0, \pmb p)
\end{equation}
and its multivariate Taylor expansion:
\begin{equation}
\label{eq:ETT}
    F_j(\pmb x_0, \pmb \theta) \approx \mathcal{P}_{F_j}^k(\delta \pmb x_0, \delta \pmb  \theta)=\overline x_{je} + \sum_{|\alpha|=1}^k \dfrac{1}{\alpha!}\partial^\alpha F_j\bigg|_{({\pmb{\overline x}}_0, {\pmb{\overline \theta}})} [\delta \pmb {x_0}, \delta \pmb  \theta]^\alpha\text{.}
\end{equation}
The coefficients of the aforementioned expansion (i.e. the $\dfrac{1}{\alpha!}\partial^\alpha F_j\bigg|_{({\pmb{\overline x}}_0, {\pmb{\overline \theta}})}$ terms) across multiple orders are here collectively termed \emph{Event Transition Tensors} (ETTs). These tensors form the foundational framework of the methodology introduced in this work, serving as the principal mathematical constructs for analyzing system dynamics at varying scales of resolution. In a system of dimension $n$, with $|\pmb \theta|$ parameters and order $k$, the dimensionality of the ETTs is given by the combinatorial sum:
$$
\sum_{i=1}^k\begin{pmatrix}
    n+|\pmb \theta|+i-1\\
    i
\end{pmatrix} = \dfrac{(n + |\pmb \theta| + i - 1)!}{k! \cdot (n + |\pmb \theta| - 1)!}
$$
This quantity grows rapidly with increasing system complexity, quickly becoming computationally intractable. Consequently, in practical applications, it becomes necessary to judiciously select only a small subset of the state vector $\pmb x$ and the parameter vector $\pmb \theta$. 
One typically must reduce the cardinality of this subset as the order of the system increases, balancing computational feasibility against the complexity of the model. 
The proposed methodology is particularly valuable when the parameters of the NeuralODE $\tilde{\pmb\theta}$ have been trained (using gradient information) and are considered as fixed, allowing for the characterization of high-order sensitivities of the NeuralODE with respect to the state $\pmb x$ and physical parameters $\pmb p$.
While event transition tensors (ETTs) can be computed considering $\tilde{\pmb\theta}$, the curse of dimensionality rapidly limits the feasibility of higher-order expansions with respect to all of the network parameters.

\begin{table}[tb]
\centering
\caption{Taxonomy of NeuralODE considered. The state dimensionality $n$, the number of parameters in the NeuralODE $|\pmb \theta|$ as well as the order used in the expansion $k$ are also reported.}
\begin{tabular}{lccccccccc}
\toprule
Problem &  \makecell{Hamiltonian\\NeuralODE} &  \makecell{Neuro\\controlled}  &\makecell{Data\\driven} & \makecell{Neural event\\manifold} & $n$ & $|\pmb \theta|$ &$k$ \\
\midrule
Lotka-Volterra   &  &  & \checkmark & &4&18&4\\
Binary asteroids      & \checkmark &  &  \checkmark & & 6 & 941 & 4\\
Comet landing    &  & \checkmark &  & \checkmark&7&2500&4\\
Drone racing   & & \checkmark &  &&16&2788&4 \\
\bottomrule
\end{tabular}
\label{tab:problems}
\end{table}

\begin{figure*}[p!]
    \centering
    \includegraphics[width=15cm]{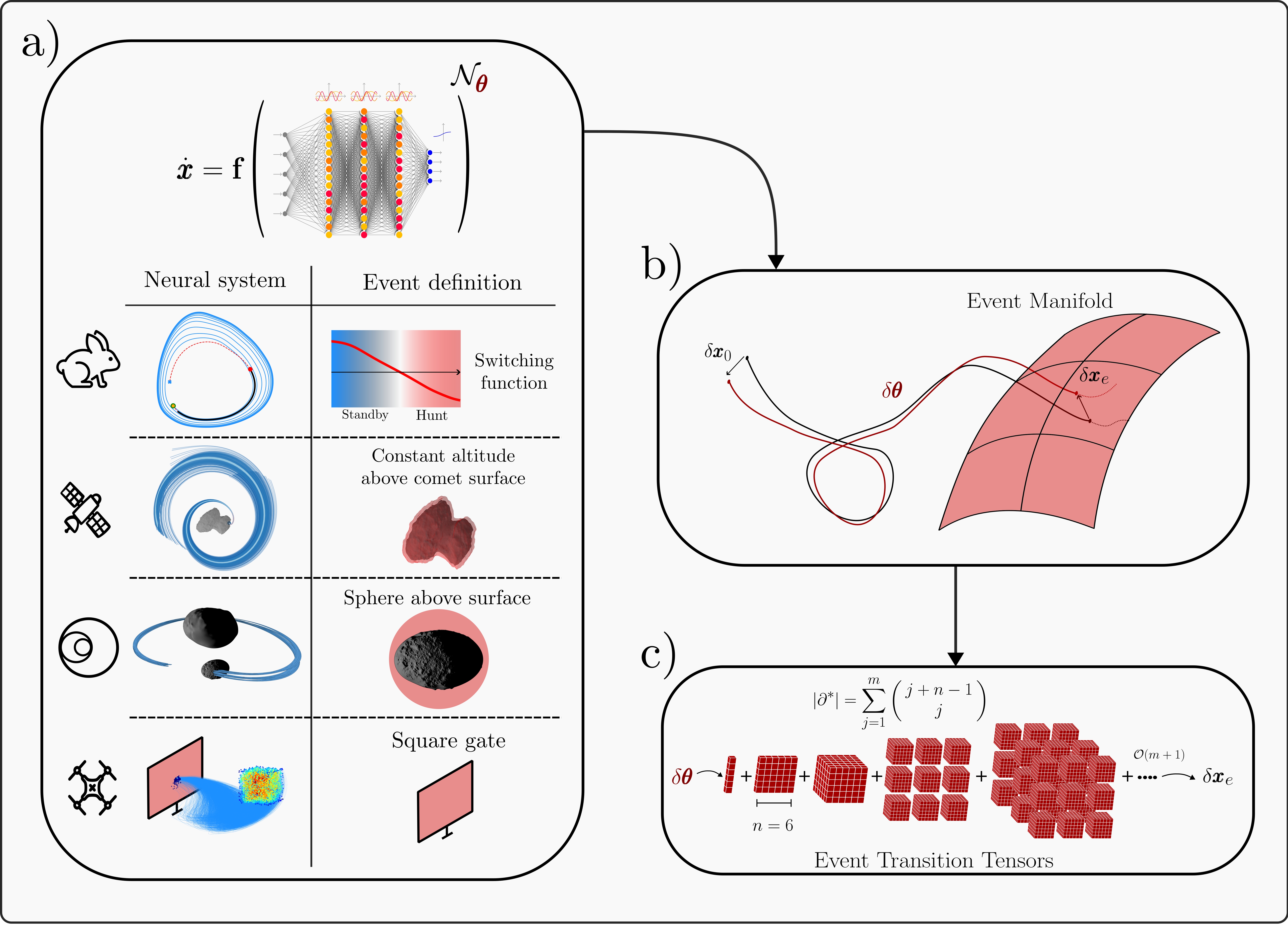}
    \caption{Schematics of the high-order expansion of a NeuralODE flow. The dynamics of a system, controlled as well as free, is learned from data resulting in a highly parametrized phase portrait which is constructed to be  differentiable everywhere by the choice of activation functions (left). The resulting highly non-linear and complex map from initial conditions to the conditions at some event manifold are represented, up to high-order, via the Event Transition Tensors (right).\label{fig:acc_err}}
\end{figure*}
\subsection*{High-order analysis of the flow}
Our methodology is outlined in Figure \ref{fig:acc_err}. 
It begins with a generic NeuralODE, an initial condition and a differentiable manifold that characterizes the triggering of some event. 
While this can include the standard case of reaching a predefined future time (in which case the ETTs collapse to the STTs and the event trigger time $t^*$ is a constant), it more interestingly covers more complex scenarios, such as reaching a subset of the system’s possible configurations.
An analytical model of the system's behavior at the event manifold with respect to variations of a subset of chosen states and parameters is thus derived by computing the corresponding ETTs defined in Eq.\eqref{eq:ETT}. The resulting expansion allows for a mathematically grounded analysis of the system's behavior, eliminating the need to rely on Monte Carlo simulations or other approximate analysis techniques. 
Depending on the system's specifics one can, for example, gain insights into the location of discontinuities of the optimal control of a NeuralODE (e.g. applying the method to the augmented system resulting from the application of Pontryagin's Maximum Principle \cite{pontryagin}).
Additionally, all the possible configurations the system reaches on the event manifold can be bounded. For example, the final landing conditions on planetary bodies of a spacecraft under uncertain initial conditions and system parameters. 
In the case of neurally controlled systems, our approach enables the verification of the stability and the identification of control margins over chosen uncertainties. 
We are also able to manipulate seamlessly generic uncertainties (i.e. also of a generic non-Gaussian type) following recent results reported in the context of astrodynamics research~\cite{acciarini2024nonlinear} which, in our case,  requires the simple evaluation of polynomial expressions involving ETTs components. 

The capabilities of our approach are best illustrated on specific systems. 
We thus have applied it to four different use cases summarized in Table~\ref{tab:problems}. The use cases are derived from two complementary paradigms that produce neuralODEs: data-driven discovery \cite{buzhardt2024relationship} and model-based optimal neural control \cite{izzo2024optimality}. 
In the realm of data-driven discovery, two different applications are explored. First, we investigate predator-prey dynamics, where we construct and analyze the (discontinuous) optimal control problem of the dynamics learned from observational data. 
Second, we address Hamiltonian learning by reconstructing governing equations for chaotic three-body gravitational dynamics from orbital observations. This case is inspired by the European Space Agency's HERA mission~\cite{michel2022esa} and proposes the use of the Hamiltonian NeuralODE framework of Greydanus et al.~(2021)~\cite{greydanus2019hamiltonian} in a space engineering context. 
In the domain of model-based optimal neural control, we focus on two advanced cases. 
The first is inspired by the problem of optimally landing a spacecraft on the weak and irregular gravitational environment of the comet 67P/Churyumov–Gerasimenko, where we also make use of an additional neural model implicitly representing the event manifold. 
In this case, we use our methodology to characterize the dispersion over the landing conditions. The second problem is drone racing, where the onboard network encodes the optimal control for a gate crossing task.

\begin{figure*}[p!]
    \centering
    \includegraphics[width=10cm]{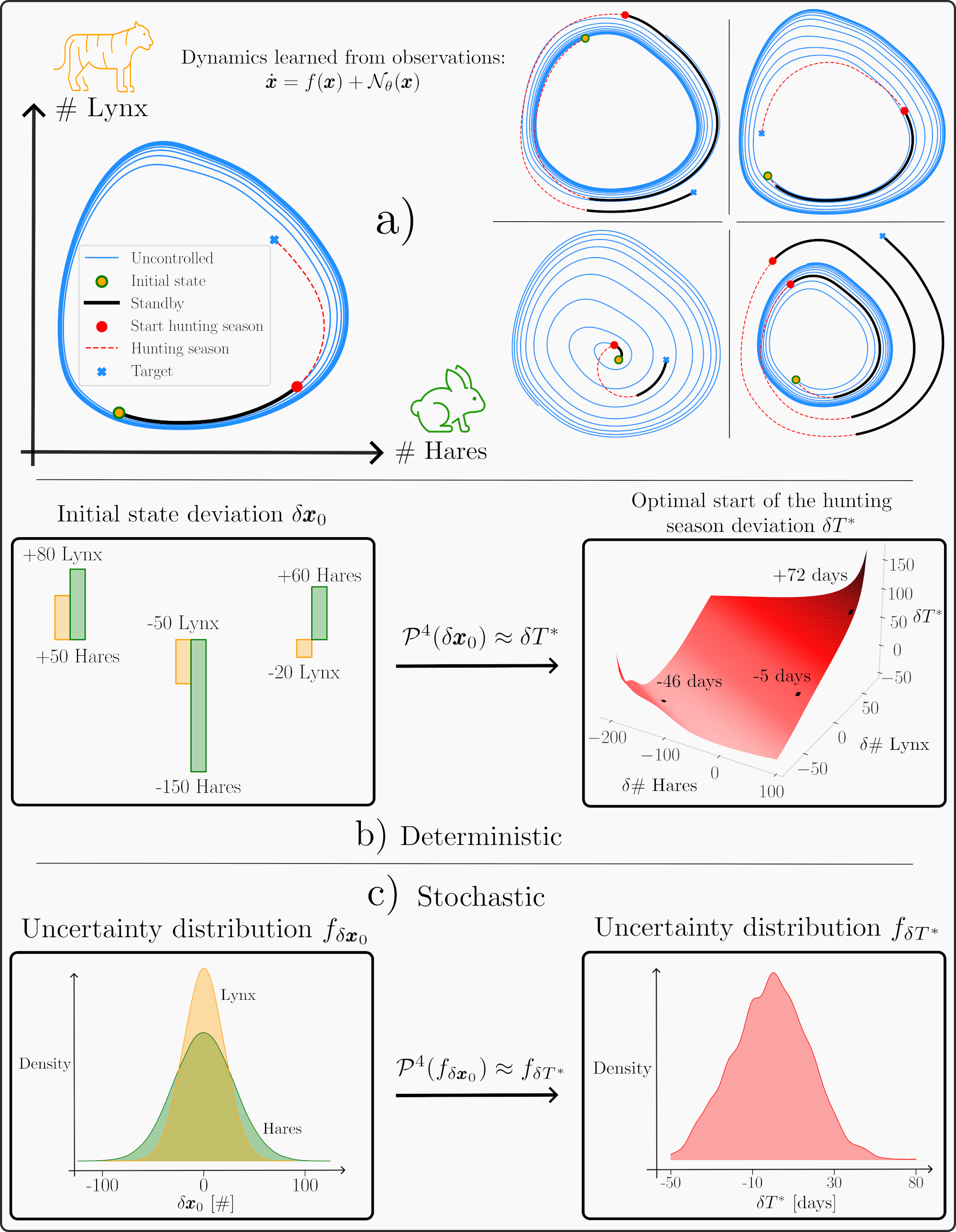}
    \caption{Learning dynamics $\dot{\pmb{x}}$ from observations in the Lotka Volterra equations. 
    a) Time evolutions of uncontrolled and controlled NeuralODEs, where the dynamics are learned from observations. Each trajectory starts at the initial state (population of Hares and Lynx) and reaches the desired target populations in a time-optimal fashion by hunting (or staying on standby) over specific time windows. On the left is the nominal trajectory around which all Taylor expansions in b) and c) are computed. On the right are examples of other time-optimal trajectories with different initial and final populations of Hares and Lynx. b) Deviations in initial populations can be mapped directly to the corresponding optimal start of the hunting season. For instance, if the initial populations of Hares and Lynx in the nominal case are off by +50 and +80, respectively, the optimal start of the hunting season is 72 days after the nominal starting date. 
    c) An uncertainty distribution in the initial state (here Gaussian) can be mapped directly to the statistical moments of the distribution of optimal start of the hunting season. Any probability distribution for which moment generating functions exists can be used to describe the uncertainty distribution in the initial state.
    \label{fig:LV}}
\end{figure*}
\subsection*{Lotka-Volterra}
In this first example, we extend the classical Lotka-Volterra equations~\cite{lotka1925elements}, widely used to model predator-prey interactions, by incorporating a neural model to account for unmodeled effects. To learn the underlying dynamics, we build on previous work~\cite{NODE4LV2021}, which demonstrated the effectiveness of NeuralODEs in accurately capturing the time-series dynamics of predator-prey populations. Accordingly, we train the system dynamics on observational data\footnote{Records of lynx and hare populations in Canada from 1900 to 1920: \url{https://jmahaffy.sdsu.edu/courses/f09/math636/lectures/lotka/qualde2.html}}. Given the limited dataset, the neural model is kept compact, consisting of only $|\pmb \theta| = 18$ parameters, to minimize the risk of overfitting.
We then formulate an optimal control problem aimed at steering the learned system from a given set of initial conditions to target conditions in minimal time, while regulating prey removal at a predefined maximum rate. The objective is to manage population dynamics (in this case, lynx and hare populations) by defining specific periods during which a maximum level of hunting is permitted. This approach provides actionable insights that can inform policy makers in developing regulations to ensure balanced and sustainable ecosystems.
To solve this control problem, we apply Pontryagin’s Maximum Principle~\cite{pontryagin} and derive a second NeuralODE representing the augmented system on which we define a two-point boundary value problem (TPBVP). This TPBVP is then solved numerically using a single-shooting method.
The solution reveals an optimal structure characterized by alternating phases: periods where the populations evolve naturally according to the learned dynamics, and controlled intervals—referred to as "hunting seasons"—during which external interventions influence the system. Further details on the learning process and the underlying mathematical framework can be found in Supplementary Note S1.

At this stage, we have defined the NeuralODE and can proceed to analyze it using our proposed methodology. Specifically, we consider a nominal trajectory generated by the optimally controlled dynamics, starting from a population of 24,700 hares and 8,600 lynx—values corresponding to the final data point in the available dataset.
As the event manifold $e$, we choose the switching function associated with the optimal control, whose zeros mark the start and end of hunting seasons. We then compute the event transition tensors up to order $k = 4$ for the corresponding event map, accounting for variations in the initial conditions. This yields an analytical model that can be used, for instance, to determine that if the initial populations of hares and lynx increase by +50 and +80, respectively, the optimal hunting season start is delayed by 72 days.
Given the inherent uncertainty in species census data, we further assume that the initial populations are uncertain. By applying exact uncertainty propagation techniques to the Taylor expansion of the event map, following the methodology described in \cite{acciarini2024nonlinear} (and briefly outlined in the Methods section), we construct the probability distribution for the optimal start of the hunting season. This effectively provides an analytical characterization of how population uncertainties impact the resulting optimal control strategy.
The results are summarized in Figure \ref{fig:LV}, where we assume that the initial population uncertainties follow normal distributions, with standard deviations of $20$ for the lynx and $30$ for the hares. Propagating these probability density functions through the TPBVP yields the resulting distribution for the optimal start of the hunting season. Notably, this distribution (see panel c of Figure \ref{fig:LV}) exhibits a skewed profile, biased toward earlier hunting start dates.

\begin{figure*}[p!]
    \centering
    \includegraphics[width=10cm]{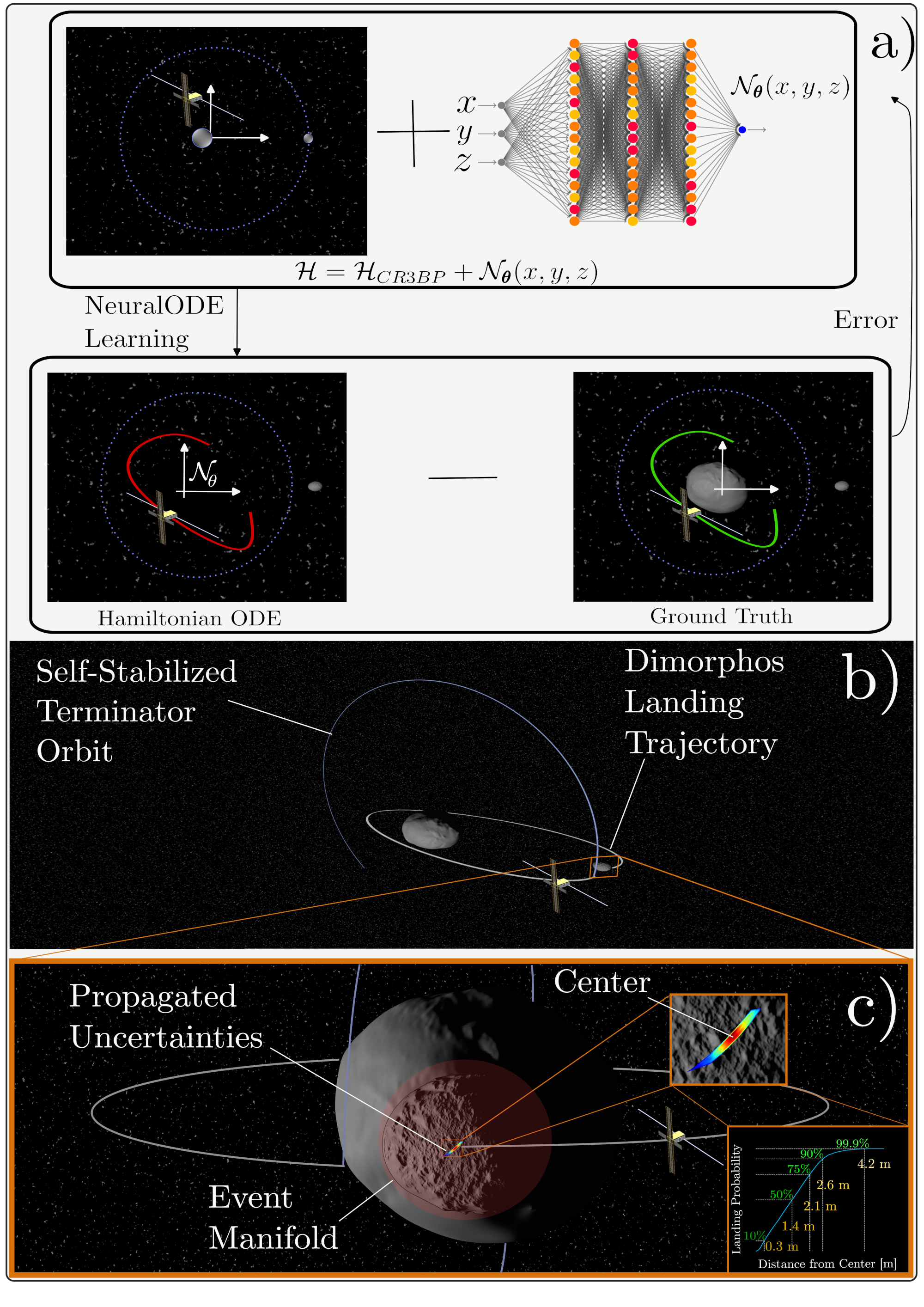}
    \caption{\textbf{Flow Expansion of Neural Hamiltonian ODE for Didymos' Irregular Gravitational Field.}  a) The upper section illustrates the learning framework, where a Neural Hamiltonian model is trained to approximate the irregular gravitational field of Didymos. The total Hamiltonian is expressed as the sum of the Circular Restricted Three-Body Problem (CR3BP) Hamiltonian (\(\mathcal{H}_{CR3BP}\)) and a neural network term (\(\mathcal{N}_{\pmb{\theta}}\)). The neural network is trained to capture deviations from the CR3BP model. The lower section represents the NeuralODE-based training setup, which backpropagates errors to refine the learned Hamiltonian.  b) The diagram presents a self-stabilizing terminator orbit (SSTO) in blue and a landing trajectory (gray) that connects the SSTO to Dimorphos’ surface. The Juventas CubeSat is also depicted.  c) A 3D visualization of the landing trajectory from b), illustrating the spatial distribution of uncertainties at the event surface, defined at a fixed altitude above Dimorphos. The color map on the landing points indicates deviations from the nominal landing position. \label{fig:image_cr3bp}}
\end{figure*}
\subsection*{Tidally locked binary asteroids}
Next, we apply our methodology to a Hamiltonian NeuralODE in the context of the Didymos-Dimorphos binary asteroid system. Binary systems make up 15\% of near-Earth objects and are key to understanding Solar System evolution \cite{richardson2006binary}. Didymos-Dimorphos has been extensively studied by NASA’s DART mission and will be the final target of ESA’s Hera mission, focusing on planetary defense \cite{rivkin2021double, michel2022esa}.

We explore a binary small body system analogous to Didymos-Dimorphos, where both bodies are assumed to be tidally locked, meaning the secondary co-rotates with the primary, maintaining synchronization with the synodic rotation period. We investigate the use of Hamiltonian ODEs to model the system’s irregular gravitational anomalies and employ a high-order expansion of the Hamiltonian ODE flow to characterize trajectory uncertainties on the secondary body’s surface.
The system extends the classical restricted circular three-body problem (CR3BP) \cite{poincare1892methodes} with an added perturbation coming from the primary’s irregular gravitational field. Assuming no prior knowledge of the primary’s shape and density, we reconstruct these irregularities from spacecraft trajectory data using a data-driven approach. Specifically, the CR3BP is modelled starting from its classical Hamiltonian and introducing an additional perturbation described by a neural network, thus forming a Hamiltonian NeuralODE \cite{hamiltonianneuralodes}. This learning process is schematically shown in Fig.~\ref{fig:image_cr3bp} a) and detailed in Supplementary Note S2.

Once the system dynamics are learned, we consider a spacecraft, modeled after Juventas CubeSat, in a self-stabilized terminator orbit (SSTO) \cite{goldberg2019juventas, scheeres2013design} and design a nominal ballistic landing trajectory from SSTO to Dimorphos’ surface, as planned for Hera \cite{goldberg2019juventas}. The resulting nominal orbits appear in Fig.~\ref{fig:image_cr3bp} b). 
Let us analyze the landing dispersion of this nominal trajectory assuming uniform (non-Gaussian) initial uncertainties of $\pm$3.5 meters in the $x$ and $y$ components of the CubeSat’s starting position.
To achieve this, a spherical event manifold, sized according to Dimorphos' semi-major axis (approximately 180 m~\cite{raducan2024physical}), is introduced, and the ETTs are computed to capture the influence of initial uncertainties on the final dispersion, as shown in Fig.\ref{fig:image_cr3bp}c. High-order uncertainty propagation~\cite{acciarini2024nonlinear} is then applied to quantify the probability of Juventas landing within a specified radius of the nominal landing trajectory. As depicted in Fig.~\ref{fig:image_cr3bp}, the model describes how  $10\%$ of trajectories land within $0.3$ m, $75\%$ within $2.1$ m, $90\%$ within $2.6$ m, and $99.9\%$ within $4.2$ m. These results provide accurate quantitative information on landing probabilities while inherently accounting for the neurally driven system dynamics.

\begin{figure*}[p!]
    \centering
    \includegraphics[width=10cm]{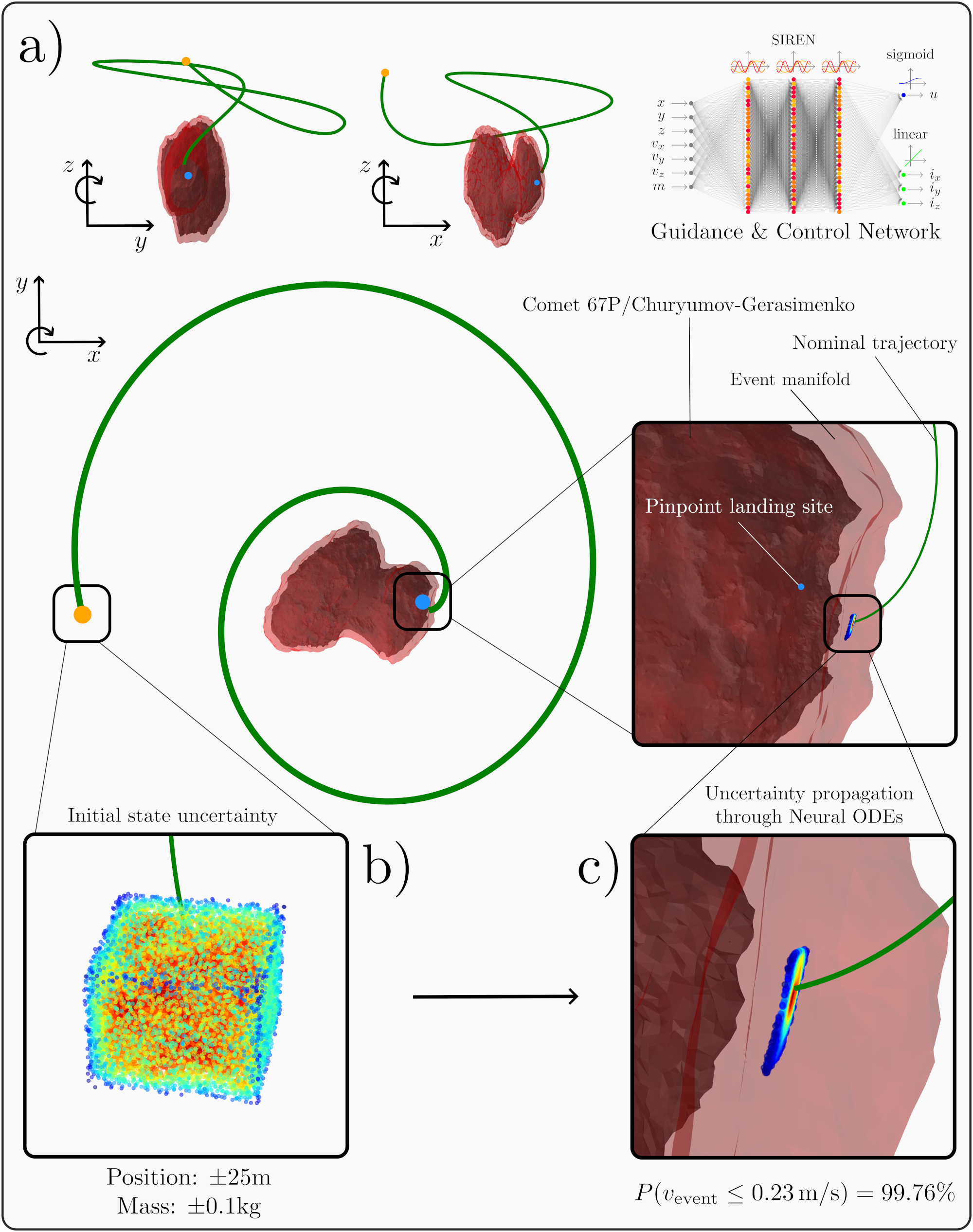}
    \caption{Mass-optimal, pinpoint spacecraft landing on the comet 67P/Churyumov-Gerasimenko using a Guidance \& Control Network (G\&CNET). The optimal control policy is learned via behavioral cloning on a dataset of optimal trajectories (not shown in this figure).
    a) The nominal neurocontrolled trajectory is shown from various perspectives in a rotating frame in which the asteroid is kept fixed. On the top right is the G\&CNET architecture, taking in as inputs the state of the spacecraft (position $[x,y,z]$, velocity $[v_x,v_y,v_z]$ and mass $m$) and outputting the corresponding optimal controls (throttle $u$ and thrust direction $[i_x,i_y,i_z]$). The G\&CNET's performance can be evaluated by computing how a given distribution in initial state uncertainty b) propagates to a specific event manifold c). 
    The event manifold used for this case is parameterized with a feedforward neural network and represents the boundary at $250$m in altitude above the surface of the comet. The color of the points on the event manifold represents the local density of points, with red corresponding to the highest density and blue to the lowest.\label{fig:67P}}
\end{figure*}
\subsection*{Comet landing}
As a third example, we address the problem of autonomous landing on a small Solar System body. Specifically, we analyze the dynamics resulting from using a neural model to represent the optimal feedback for a mass-optimal pinpoint landing on the highly irregular surface of comet 67P/Churyumov-Gerasimenko, visited in 2014 by the European Space Agency's Rosetta mission~\cite{glassmeier2007rosetta}.

Such models, called Guidance \& Control Networks (G\&CNETs), have gained attention as a viable alternative to Model Predictive Control (MPC) methods in the field of space research and drone racing \cite{DNN_improved_tracking, osti_10100589, Kaufmann2023, aggressive_online_control, cheng2018real, federici2021deep, dario_seb_gcnet, origer2023quad, origer2024neuralODEfix, izzo2020stability, izzo2021real}. 
By directly translating system states into optimal control inputs, they integrate guidance and control into a single inference step. 
The two primary learning paradigms for training G\&CNETs are reinforcement learning and behavioral cloning~\cite{izzo2024optimality}.
While in this work we make use of behavioral cloning, thus leveraging a dataset of optimal trajectories, our methodology is agnostic to how optimality principles are embedded into the neural network.

To construct the dataset of optimal trajectories, we first generate a nominal trajectory representing the spacecraft's designed landing phase by solving the TPBVP derived from Pontryagin’s Maximum Principle~\cite{pontryagin}. Next, we extend the dynamics to the augmented system and efficiently generate millions of optimal trajectories using a technique known as backward generation of optimal examples (BGOE)~\cite{izzo2021real, dario_seb_gcnet, origer2024neuralODEfix}.
The resulting dataset consists of optimal trajectories that originate from different initial conditions but converge at the desired pinpoint landing site. We finally train a G\&CNET to learn the optimal throttle input and thrust direction required for precision landing.
Further details on the dataset generation and the training process are provided in Supplementary Note S3.

The neural model employed for this case is a SIREN \cite{sitzmann2020implicit} G\&CNET parameterized with $|\pmb \theta| = 2,500$ parameters, which allows high accuracy when learning a representation for the optimal feedback \cite{origer2024gcnetperiodic}. 
Accuracy is important here, whereby the G\&CNET, acting as a state-feedback, needs to cancel out prediction errors during the landing phase faster than it accumulates them.
Once trained, it is essential to certify and validate our neuro-controlled system. 
Consider, for instance, the need to verify the following requirement:
\begin{quote}
The G\&C algorithm shall steer the spacecraft to an altitude of $250$ m $\pm 5$ m above the surface of the comet while ensuring a relative velocity of $\leq 0.5$ m/s. 
The algorithm is required to achieve this relative velocity in at least 99\% of cases, given an initial state uncertainty of $\pm 25$m in position and $\pm 0.1$kg in spacecraft mass.
\end{quote}
To analyze the system state at an altitude of $250$m above the comet's surface, we define this boundary as the event manifold. We take into account the highly irregular shape of 67P/Churyumov-Gerasimenko by representing it as an event manifold $e(\pmb x) = 0$ training a second additional neural network, (also a SIREN), which further demonstrate the versatility of our methodology in handling complex events. We thus obtain a NeuralODE, a nominal trajectory and a differentiable event manifold and can compute the corresponding ETTs. The Taylor expansion of the event map allows us to analyze the performance of the neurally controlled dynamics, as illustrated in Figure \ref{fig:67P}. In particular, we can introduce non-Gaussian uncertainties on the initial conditions defining a uniform probability density function over the initial state of the magnitude set in the requirement.
Propagating this uncertainty at the event, it is then possible to compute that given an initial state uniform uncertainty of $\pm25$m in position and $\pm0.1$kg in mass, the final relative velocity at the event (altitude of $250$ meters from the comet) is less than $0.23$m/s in $99.76\%$ of cases. Hence certifying that the neurally controlled dynamics meets the engineering requirement defined above. 

\begin{figure*}[p!]
    \centering
    \includegraphics[width=10cm]{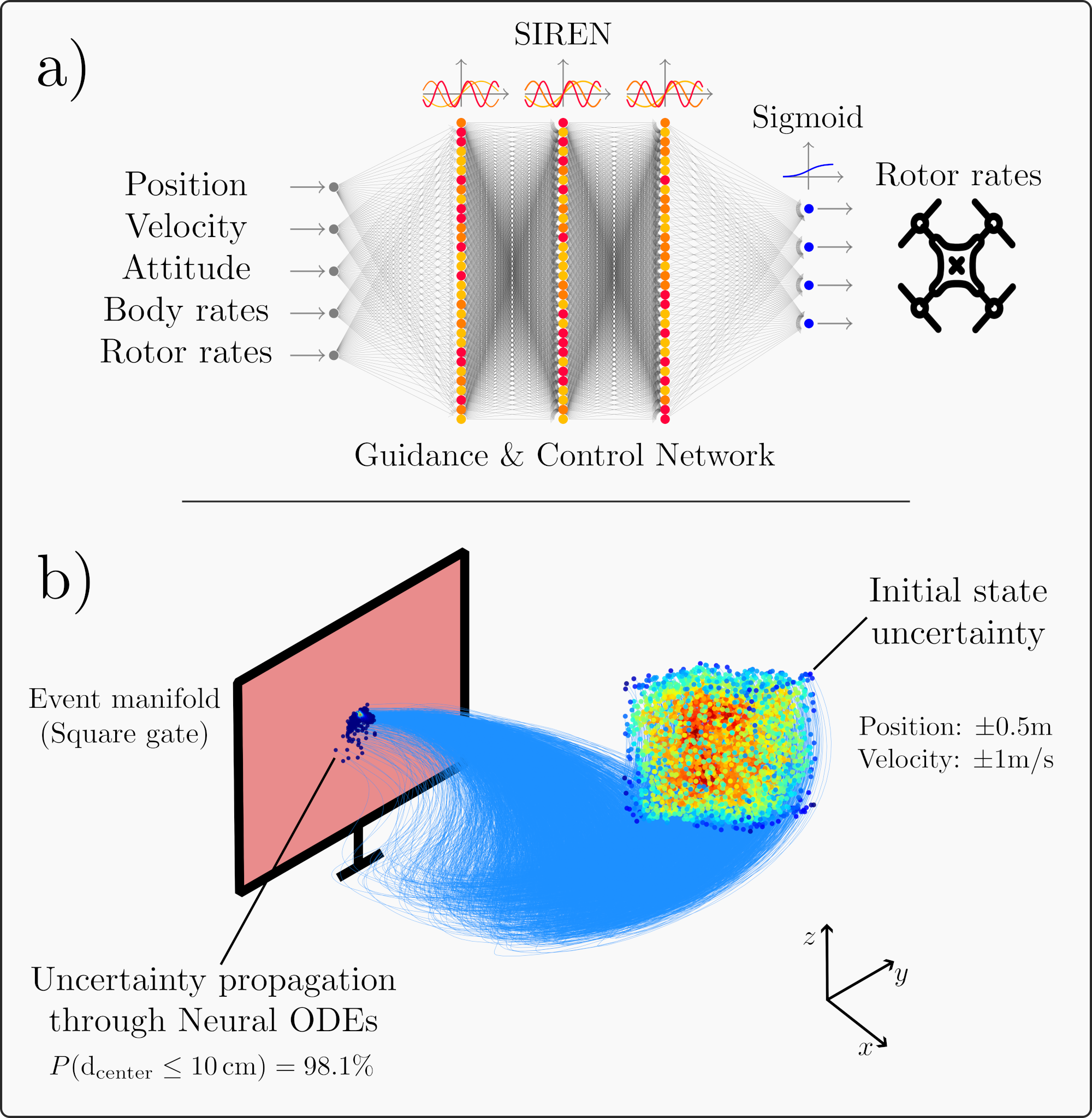}
    \caption{Energy-optimal drone racing through a square gate using a Guidance \& Control Network (G\&CNET). The optimal control policy is learned via behavioral cloning on a dataset of optimal trajectories (not shown in this figure).
    a) The G\&CNET architecture, taking in as inputs the full state vector of the drone and outputting the corresponding optimal controls. 
    b) The G\&CNET's performance can be evaluated by computing how a given distribution in initial state uncertainty  propagates to a specific event manifold. For instance, it can verified that, given an initial state uncertainty of $\pm0.5$m in position and $\pm1$m/s in velocity, the final distance to the center of the gate at the event, $d_{\text{center}}$, is less than $10$cm in $98.1\%$ of cases. The event manifold used in this case is a simple two-dimensional square gate. The color of the points on the event manifold represents the local density of points, red corresponding to the highest density and blue to the lowest.
    \label{fig:drone_racing}}
\end{figure*}
\subsection*{Drone racing}
In this last example, we obtain a NeuralODE by training a G\&CNET to guide and control a quadcopter energy-optimally through the center of a two-dimensional square gate. 
Quadcopters pose a particularly challenging testbed for G\&CNETs due to their high-dimensional state space (16 states) and the nonlinear dynamics arising from aerodynamic effects. Notably, drone racing is used as a gym environment to evaluate G\&CNETs on physical robotic platforms, thereby increasing confidence in their potential application to future space exploration missions \cite{izzo2024optimality}. 
In contrast to the mass-optimal trajectory generation for comet landings, where indirect methods are employed, here we leverage direct optimization to construct a dataset of optimal quadcopter trajectories \cite{origer2023quad}. Importantly, our approach remains agnostic to the specific training framework used for the network, so that the same applies if reinforcement learning—a commonly favored strategy in drone racing \cite{Kaufmann2023}—was used to train the G\&CNET. More details on the dataset generation and learning process are provided in Supplementary Note S4.
Similarly to the comet landing scenario, the neural model employed is a SIREN \cite{sitzmann2020implicit} G\&CNET containing $|\pmb \theta| = 2,788$ parameters, allowing a high accuracy in the optimal feedback representation \cite{origer2024gcnetperiodic}.
G\&CNETs are well suited for problems with fast changing dynamics as the low computational power associated with a single inference allows to run them at high frequency. 
For reference, G\&CNETs of similar size to ours have been tested on the Parrot Bebop 1 which uses a Parrot P7 dual-core Cortex A9 CPU allowing it to infer the network at a frequency of $450$Hz \cite{origer2023quad}.

Following the training of the G\&CNET, we obtain a closed-loop dynamical system describing the quadcopter's motion, where the control inputs (rotor rates) are generated by the G\&CNET. This yields a NeuralODE formally expressed in the form of Eq. ($\ref{eq:ETT}$). To analyze the system's neural-driven dynamics, we define the event manifold as the plane coinciding with the two-dimensional gate (see Figure \ref{fig:drone_racing}) and subsequently compute the associated event transition tensors. 
To evaluate the system's robustness under initial state perturbations, we introduce a uniform uncertainty characterized by position deviations of $\pm 0.5$ m and velocity deviations of $\pm 1$ m/s. Propagating this uncertainty through the neural flow yields an analytical representation of the probability density function governing the gate-crossing event. This can be leveraged to verify key performance metrics. For example, we verify that the distance from the quadcopter to the gate center at the crossing event is less than 10 cm with a probability of 0.981.

\section*{Discussion}
Our methodology analyses the performance of NeuralODEs by constructing a high-order Taylor representation of their behavior at the time, $t^*$ which may correspond to either a fixed time instant or, more significantly, the variable time at which a specific event occurs. The practical applicability of the approach depends on both the feasibility of computing high-order tensors for the given NeuralODE and the accuracy of the resulting Taylor expansion.

Regarding computational feasibility, in the four cases examined in this study, careful selection of state and parameter subsets ensures that computations remain tractable, typically completing within minutes while providing valuable insights into system dynamics. However, increasing problem complexity—such as extending to higher-order representations or incorporating additional uncertainties—would significantly increase dimensionalities and thus require more computational resources.

With respect to accuracy, previous studies employing high-order expansions of ordinary differential equation flows \cite{berz1988differential, wittig2015propagation, acciarini2024nonlinear, boone2021orbital, izzo2020stability} have investigated techniques for estimating the remainder of the multivariate Taylor expansions and discussed the issue of their rigorous validation. Numerical methodologies have been proposed and developed to ensure a prescribed level of accuracy. 
It is worth noting that in this study, the accuracy of the Taylor expansions proved sufficient to obtain the results reported, as confirmed through Monte Carlo simulations of the full nonlinear NeuralODE dynamics.
Nevertheless, a fully rigorous mathematical proof of the certification and accuracy of the computations remains, in the general case, an open challenge.
Techniques such as Taylor models \cite{berz1998verified} offer a promising avenue to achieve such a formal verification, but also introduce additional computational complexity.

Our method marks a significant advancement in the analysis of NeuralODEs by providing an analytical representation of their behavior on differentiable event manifolds within a well-defined region around a nominal trajectory. 
This representation contributes to a rigorous and computationally efficient alternative to extensive Monte Carlo simulations, allowing for a precise assessment of system behavior in the presence of unmodeled effects or uncertainties. In safety-critical applications, such as aerospace and automotive systems, this approach is expected to enhance trust in neurally driven dynamics and contribute to the rigorous verification and certification of their reliability.

\section*{Methods}

\subsection*{Computing ETTs}
At the core of our findings lies the efficient computation of ETTs, as defined in Eq. \eqref{eq:ETT}. 
Throughout this work, we leverage a state-of-the-art implementation of arbitrary-order variational equations within a Taylor integration framework \cite{biscani2022reliable, biscani2021revisiting}, made openly available in the \emph{heyoka} software package~\cite{biscani2025heyoka}. 
While high-order variational equations have been criticized for their derivational complexity \cite{park2006nonlinear} and the increased system complexity compared to alternative algebraic approaches, such as forward-mode automatic differentiation applied directly to ODE solvers \cite{valli2013nonlinear}, our results demonstrate their effectiveness and practical feasibility. 
The derivation of the variational equations for differentiable systems (hence also for all our NeuralODEs) can be fully automated by leveraging efficient symbolic manipulation of expressions. 
Additionally, the integration of just-in-time compilation (via LLVM) \cite{lattner2004llvm} and careful sub-expression reuse, combined with an underlying modern Taylor integration scheme, makes this technique particularly efficient in our context.

Consider a general NeuralODE of the form:
$$
\left\{
\begin{array}{l}
\dot{\pmb x} = \pmb f(\pmb x, \pmb \theta, t) \\
\pmb x(0) = \pmb x_0
\end{array}
\right.
$$
where the system evolves until it reaches a differentiable event manifold, defined by $e(\pmb x, \pmb  \theta, t) = 0$, at time $T$. To recast this as a fixed-time integration problem, the normalized time variable $\tau = t / T$ is introduced. An additional equation to track the time variation $\epsilon = \frac {de}{dt}$ of the event function  is also added leading to the an augmented system of ODEs:

$$
\left\{
\begin{array}{l}
\dot{\pmb x} = T \pmb f(\pmb x, \pmb  \theta, \tau) \\
\dot \varepsilon = T\dfrac{\partial e}{\partial t} + T \nabla_{\pmb x} e \cdot f \\
\pmb x(0) = \pmb x_0 \\
\varepsilon(0) = e(\pmb x(0), \pmb  \theta, 0)
\end{array}
\right.
\text{.}
$$
By integrating the order $k$ variational equations associated to the above system up to $\tau=1$, we obtain the necessary elements to construct a polynomial representation of the system state at the final integration time:

$$
\left\{
\begin{array}{l}
\delta \pmb x_f = \mathcal P^k_{\pmb x_f}(\delta \pmb x_0, \delta \pmb  \theta, \delta T) \\
\delta \varepsilon_f = \mathcal P^k_{\varepsilon_f}(\delta \pmb x_0, \delta \pmb  \theta, \delta T, \delta \varepsilon_0)
\end{array}
\right.
\text{,}
$$
where $\mathcal P^k$ denotes generic multivariate polynomials of order $k$. If the implicit function theorem \cite{krantz2002implicit} holds, this polynomial system can be efficiently and exactly inverted. Using the exact map inversion algorithm \cite{berz1998verified}, we can swap variables between the left- and right-hand sides. Specifically, swapping $\varepsilon_f$ and $\delta T$ in the last equation and forcing the final event variation to zero leads to write:
$$
\delta T^* = \mathcal P^k_{\varepsilon_f}(\delta \pmb x_0, \delta \pmb  \theta, \cancelto{0}{\delta \varepsilon_f}, \delta \varepsilon_0)
\text{,}
$$
which expresses the variation in event trigger time as a function of perturbations in the initial conditions and system parameters. Substituting this result into $\mathcal P^k_{\pmb x_f}$, we obtain an order-$n$ representation of the system state constrained to the event manifold, precisely defining the ETTs.
Manually performing these derivations and computations for each system would be impractical. To address this, we have developed a fully automated pipeline that seamlessly applies the above methodology to any generic system of ODEs, including NeuralODEs. This implementation is available within the codebase of the open-source \emph{heyoka} project.

\subsection*{Propagation of uncertainties through NeuralODEs}
As shown in Eq.~\eqref{eq:ETT}, the computation of ETTs enables us to construct a polynomial analytical model that describes the system’s state at the event manifold. This model captures how changes in the system parameters and initial conditions relate to variations in the final state at the event. Due to its polynomial nature, it allows for efficient uncertainty propagation. Specifically, we leverage the extension of Park \& Scheeres' work to non-Gaussian uncertainties~\cite{acciarini2024nonlinear, park2006nonlinear} to propagate uncertainties within our NeuralODE framework. The underlying idea can be summarized as follows: by expressing the state evolution through ETTs, we encode the sensitivity of the system's dynamics to initial uncertainties and parameter variations in a compact, higher-order form. This enables accurate and computationally efficient propagation of uncertainties, particularly in nonlinear regimes where traditional first-order methods would fall short.
While the notation and formalism in high-dimensional, multivariate cases can become cumbersome, here we illustrate the procedure using a simple one-dimensional system with a single parameter. 
For an extensive treatment of the uncertainty propagation method here used we refer to previous work \cite{acciarini2024nonlinear}.
We assume that the ETTs have been computed up to the second order. Leveraging Eq.~\eqref{eq:ETT}, the resulting truncated Taylor expansion of the event map can thus be written as:
$$
\delta x_e = \frac{\partial x_e}{\partial x_0} \delta x_0 + \frac{\partial x_e}{\partial \theta} \delta \theta + \frac{1}{2} \frac{\partial^2 x_e}{\partial x_0^2} (\delta x_0)^2 + \frac{1}{2} \frac{\partial^2 x_e}{\partial \theta^2} (\delta \theta)^2 + \frac{\partial^2 x_e}{\partial x_e \partial \theta} \delta x_0 \delta \theta,
$$
Here, $x_0$ and $x_e$ represent the initial state and the state at the event, $\theta$ is the system parameter, and the partial derivatives are the various ETT components.
We can now propagate statistical moments, which characterize uncertainties in the initial state and parameters, directly onto the event manifold by evaluating them over the polynomial described above. For instance, in the case of a second-order moment, we have:
$$
m_2 = \mathbb E[(\delta x_e-\mathbb E(\delta x_e))(\delta x_e-\mathbb E(\delta x_e))] = \mathbb E(\delta x_e^2) - \big(\mathbb E(\delta x_e)\big)^2
$$
where $\delta x_e^2$ is a fourth order polynomial in $\delta x_0, \delta \theta$, and we find:
$$
m_2 = \sum_{|\alpha+\beta| = 2 }^4 a_{\alpha \beta} \mathbb E[\delta x_0^\alpha \delta \theta^\beta] - \big(\mathbb E(\delta x_e)\big)^2
$$
The coefficients $a_{\alpha \beta}$ are computed via polynomial multiplication and solely depend on the ETTs. The various moments of the uncertainties distribution  $\mathbb E[\delta x_0^\alpha \delta \theta^\beta]$ are efficiently precomputed using symbolic manipulation of the corresponding moment-generating function. Notably, these moments depend only on the chosen distribution, not on the specific NeuralODE under consideration. The remaining term $\big(\mathbb E(\delta x_e)\big)^2$ is linked to the computation of the lower-order moment $m_1$. The same approach can be used for higher-dimensional systems and higher-order moments. The key takeaway is that ETTs enable the efficient determination of all moments of the distribution of the NeuralODE state at the event manifold.

\subsection*{Neural architectures}
We employ the SIREN architecture by Sitzmann et al. \cite{sitzmann2020implicit} for all our neural networks in this work.
SIREN networks provide several benefits, including being differentiable everywhere and displaying remarkable approximation power on a wide range of regression tasks ranging from image reconstruction to guidance \& control and implicit shape representations \cite{sitzmann2020implicit, izzo2022geodesy,origer2024gcnetperiodic}.

The architecture consists of a simple feed-forward neural network with sine functions as activation for the hidden layers:
\begin{equation}\label{eq:SIREN}
\text{SIREN}: \left\{ 
\begin{array}{l}
        \Phi(\mathbf{x}) = \mathbf{W}_n(\phi_{n-1}\circ\phi_{n-2}\circ \dots \circ \phi_0)(\mathbf{x})+\mathbf{b}_n \\
        \mathbf{x}_i \mapsto \phi_i(\mathbf{x}_i)=\sin{(\omega_0\mathbf{W}_i\mathbf{x}_i+\mathbf{b}_i)}
\end{array}
\right.
\end{equation}
The $i^{th}$ layer of the neural network is represented by the function $\phi_i:\mathbb{R}^{M_i} \to \mathbb{R}^{N_i}$, where $\mathbf{W}_i \in \mathbb{R}^{N_i \times M_i}$ is the weight matrix, $\mathbf{b}_i \in \mathbb{R}^{N_i}$ is the bias vector, and $\mathbf{x}_i \in \mathbb{R}^{M_i}$ is the input vector.
Non-linearity is introduced by applying the sine function to each component of the resulting vector.
We follow the suggested weight initialization scheme and input normalization laid out in \cite{sitzmann2020implicit}.
Specifically, we use $\omega_0=30$ and initialize the weights using $w_i\sim \mathcal{U}(-1/n,1/n)$ for the first layer and $w_i\sim \mathcal{U}(-\sqrt{6/n}/\omega_0,\sqrt{6/n}/\omega_0)$ for all other layers, where $n$ is here the number of inputs to the $i$-th neuron. Furthermore, we add an input normalization layer such that the minimum and maximum value for each input corresponds to $-1$ and $1$, respectively. The use of SIREN networks enables the neural models to remain relatively small in size \cite{origer2024gcnetperiodic}, which is crucial for enhancing the computational efficiency of both the ETT calculations and the NeuralODE learning process. As demonstrated in Supplementary Table S1, our computational pipeline exhibits notable performance gains, particularly at lower model dimensionalities, achieving significant speedup.

\subsection*{Data-driven training}
We train the parameters of our (data-driven) NeuralODEs (see Table (\ref{tab:problems})) to reproduce observed data $(t_i, \hat{\pmb x}_i)$ by using the loss function $\mathcal L = \frac1N\sum_{i=0}^{N-1} \mathcal L_i = \frac1N\sum_{i=0}^{N-1}||\pmb \Phi(t_{i+1}; \pmb x_i,\pmb \theta) - \hat{\pmb x}_{i+1}||^2$ in a stochastic gradient descent algorithm. In essence, our objective is to minimize the prediction error arising from propagating the NeuralODE between consecutive observations. Efficient gradient computation with respect to the system parameters $\pmb\theta$ is essential for this optimization process. Various techniques have been developed to tackle this challenge across diverse contexts. Since only first-order information (i.e., the gradient) is required in this case, a widely adopted approach employs reverse-mode automatic differentiation (AD) through ODE solvers, typically based on Runge-Kutta methods \cite{torchdyn, torchdiffeq}. However, this approach faces substantial memory limitations when the number of integration steps increases, as the complexity of the computational graph grows rapidly. Moreover, it requires re-computation of the computational graph for each integration, due to dynamic variations in both the number and positioning of the time steps, further amplifying the computational overhead. 
To mitigate these limitations and implement a method agnostic to the system dynamics, data characteristics, and time grid, we instead adopt the adjoint method, as originally proposed by Chen et al. \cite{chen2018neural, torchdiffeq} and also suggested as a baseline in the associated \emph{torcheqdiff} open-source project~\cite{torchdiffeq}. Unlike that baseline, however, we leverage a modern Taylor integration scheme, as opposed to Runge-Kutta methods, to achieve additional numerical stability and computational speed-ups, following the methodology outlined in Biscani and Izzo \cite{biscani2021revisiting}. This combination enhances both efficiency and flexibility, making it particularly well-suited for NeuralODE applications requiring efficient gradient computation over varying time grids and complex dynamical systems. 

For each term in the loss function, we introduce a vector quantity $\pmb a_i$, referred to as the adjoint, defined as:
$$
\pmb a_i = \frac{\partial \mathcal L_i}{\partial \pmb x(t)},
$$
along with auxiliary vector quantities representing the loss sensitivities:
$$
\pmb \lambda_i = \frac{\partial \mathcal L_i}{\partial \pmb \theta},
$$
To compute the gradient, we consider the following augmented system of ordinary differential equations (ODEs):
\begin{equation}
\label{eq:aug_adjoint}
\left\{
\begin{array}{l}
\dot{\pmb x} = \pmb f(\pmb x, \pmb p, t) \\
\dot {\pmb a_i} = - \pmb a_i^T\nabla_{\pmb x} \pmb f\\
\dot {\pmb \lambda_i} = - \pmb a_i^T \nabla_{\pmb \theta} \pmb f
\end{array}
\right.
\end{equation}
which is integrated backward in time from the final conditions $\pmb x=\pmb x_f$, $\pmb a_i = \pmb a_{i_f}$ and $\pmb \lambda_i = \pmb 0$. The terminal state $\pmb x_f$ is obtained via an additional (forward) numerical integration of the original dynamics, starting from $\pmb x_i$ and advancing up to $t_{i+1}$. In our specific case, the terminal adjoint is given by:
$$
\pmb a_{i_f} = 2 (\pmb x_f - \hat{\pmb x}_{i+1})
$$
where $\hat{\pmb x}_{i+1}$ denotes the target state at time $t_{i+1}.$
In Supplementary Note S5, we present a quantitative benchmark of our adjoint method implementation, comparing it against the default method in \emph{torchdiffeq}. The results demonstrate that our adaptive Taylor-based approach can offer a significant advantage, particularly in systems with a moderate number of parameters, as well as in scenarios requiring high precision for gradient computation.

\subsection*{G\&CNET training}

In this work, all G\&CNETs are feed-forward neural networks trained via behavioral cloning. The training process begins with generating a dataset of optimal trajectories, which is then split into features (system states) and labels (optimal controls).  For the drone racing problem, each trajectory is solved individually using a direct method \cite{origer2023quad}. In contrast, for the mass-optimal landing, we employ the Backward Generation of Optimal Examples (BGOE) technique \cite{dario_seb_gcnet, origer2024neuralODEfix, izzo2021real}, which significantly accelerates dataset generation for problems where Pontryagin's Maximum Principle is applied \cite{izzo2024optimality}.

\section*{Code availability}
All code used to produce the results is available online at 
\url{http://bit.ly/3FPRTUh} and is licensed under an Apache 2.0 open source license.

\bibliography{scibib}

\begin{thebibliography}{10}

\bibitem{tsoulos2006solving}
I.~G. Tsoulos, I.~E. Lagaris, Solving differential equations with genetic programming.
\newblock {\it Genetic Programming and Evolvable Machines\/} {\bf 7}, 33--54 (2006).

\bibitem{schmidt2009distilling}
M.~Schmidt, H.~Lipson, Distilling free-form natural laws from experimental data.
\newblock {\it science\/} {\bf 324}, 81--85 (2009).

\bibitem{izzo2017differentiable}
D.~Izzo, F.~Biscani, A.~Mereta, {\it Genetic Programming: 20th European Conference, EuroGP 2017, Amsterdam, The Netherlands, April 19-21, 2017, Proceedings 20\/} (Springer, 2017), pp. 35--51.

\bibitem{lejarza2022data}
F.~Lejarza, M.~Baldea, Data-driven discovery of the governing equations of dynamical systems via moving horizon optimization.
\newblock {\it Scientific Reports\/} {\bf 12}, 11836 (2022).

\bibitem{greydanus2019hamiltonian}
S.~Greydanus, M.~Dzamba, J.~Yosinski, Hamiltonian neural networks.
\newblock {\it Advances in neural information processing systems\/} {\bf 32} (2019).

\bibitem{chen2018neural}
R.~T. Chen, Y.~Rubanova, J.~Bettencourt, D.~K. Duvenaud, Neural ordinary differential equations.
\newblock {\it Advances in neural information processing systems\/} {\bf 31} (2018).

\bibitem{izzo2024optimality}
D.~Izzo, E.~Blazquez, R.~Ferede, S.~Origer, C.~De~Wagter, G.~C. de~Croon, Optimality principles in spacecraft neural guidance and control.
\newblock {\it Science Robotics\/} {\bf 9}, eadi6421 (2024).

\bibitem{buzhardt2024relationship}
J.~Buzhardt, C.~R. Constante-Amores, M.~D. Graham, On the relationship between koopman operator approximations and neural ordinary differential equations for data-driven time-evolution predictions.
\newblock {\it arXiv preprint arXiv:2411.12940\/}  (2024).

\bibitem{hornik1989multilayer}
K.~Hornik, M.~Stinchcombe, H.~White, Multilayer feedforward networks are universal approximators.
\newblock {\it Neural networks\/} {\bf 2}, 359--366 (1989).

\bibitem{calin2020deep}
O.~Calin, {\it Deep learning architectures\/} (Springer, 2020, pp. 251-284).

\bibitem{bonnaffe2021neural}
W.~Bonnaff{\'e}, B.~C. Sheldon, T.~Coulson, Neural ordinary differential equations for ecological and evolutionary time-series analysis.
\newblock {\it Methods in Ecology and Evolution\/} {\bf 12}, 1301--1315 (2021).

\bibitem{bottcher2022ai}
L.~B{\"o}ttcher, N.~Antulov-Fantulin, T.~Asikis, Ai pontryagin or how artificial neural networks learn to control dynamical systems.
\newblock {\it Nature communications\/} {\bf 13}, 333 (2022).

\bibitem{kumar2023physics}
T.~Kumar, A.~Kumar, P.~Pal, {\it NeurIPS Machine Learning and the Physical Sciences Workshop\/} (2023).

\bibitem{giang2024conditional}
K.~T. Giang, Y.~Kim, A.~Finazzi, Conditional latent odes for motion prediction in autonomous driving.
\newblock {\it arXiv preprint arXiv:2405.19183\/}  (2024).

\bibitem{chen2022forecasting}
X.~Chen, F.~A. Araujo, M.~Riou, J.~Torrejon, D.~Ravelosona, W.~Kang, W.~Zhao, J.~Grollier, D.~Querlioz, Forecasting the outcome of spintronic experiments with neural ordinary differential equations.
\newblock {\it Nature communications\/} {\bf 13}, 1016 (2022).

\bibitem{griewank1994computational}
A.~Griewank, {\it Proceedings\/} (1994), vol.~1, pp. 37--53.

\bibitem{beutler2005variational}
G.~Beutler, Variational equations.
\newblock {\it Methods of Celestial Mechanics: Volume I: Physical, Mathematical, and Numerical Principles\/} pp. 175--207 (2005).

\bibitem{skokos2010numerical}
C.~Skokos, E.~Gerlach, Numerical integration of variational equations.
\newblock {\it Physical Review E—Statistical, Nonlinear, and Soft Matter Physics\/} {\bf 82}, 036704 (2010).

\bibitem{berz1988differential}
M.~Berz, {\it AIP Conference Proceedings\/} (American Institute of Physics, 1988), vol. 177, pp. 275--300.

\bibitem{corliss1997high}
G.~F. Corliss, A.~Griewank, P.~Henneberger, G.~Kirlinger, F.~A. Potra, H.~J. Stetter, {\it Numerical Analysis and Its Applications: First International Workshop, WNAA'96 Rousse, Bulgaria, June 24--26, 1996 Proceedings 1\/} (Springer, 1997), pp. 114--125.

\bibitem{griewank2008evaluating}
A.~Griewank, A.~Walther, {\it Evaluating derivatives: principles and techniques of algorithmic differentiation\/} (SIAM, 2008).

\bibitem{jax2018github}
J.~Bradbury, R.~Frostig, P.~Hawkins, M.~J. Johnson, C.~Leary, D.~Maclaurin, G.~Necula, A.~Paszke, J.~VanderPlas, S.~Wanderman-Milne, Q.~Zhang, {JAX}: composable transformations of {P}ython+{N}um{P}y programs (2018). Version 0.3.13.

\bibitem{biscani2021revisiting}
F.~Biscani, D.~Izzo, {Revisiting high-order Taylor methods for astrodynamics and celestial mechanics}.
\newblock {\it Monthly Notices of the Royal Astronomical Society\/} {\bf 504}, 2614-2628 (2021).

\bibitem{park2006nonlinear}
R.~S. Park, D.~J. Scheeres, Nonlinear mapping of gaussian statistics: theory and applications to spacecraft trajectory design.
\newblock {\it Journal of guidance, Control, and Dynamics\/} {\bf 29}, 1367--1375 (2006).

\bibitem{bani2019exact}
A.~Bani~Younes, Exact computation of high-order state transition tensors for perturbed orbital motion.
\newblock {\it Journal of Guidance, Control, and Dynamics\/} {\bf 42}, 1365--1371 (2019).

\bibitem{boone2021orbital}
S.~Boone, J.~McMahon, Orbital guidance using higher-order state transition tensors.
\newblock {\it Journal of Guidance, Control, and Dynamics\/} {\bf 44}, 493--504 (2021).

\bibitem{izzo2020stability}
D.~Izzo, D.~Tailor, T.~Vasileiou, On the stability analysis of deep neural network representations of an optimal state feedback.
\newblock {\it IEEE Transactions on Aerospace and Electronic Systems\/} {\bf 57}, 145--154 (2020).

\bibitem{berz2002cosy}
M.~Berz, K.~Makino, Cosy infinity version 8.1. user’s guide and reference manual.
\newblock {\it Department of Physics and Astronomy MSUHEP-20704, Michigan State University\/}  (2002).

\bibitem{izzo2020dcgp}
D.~Izzo, F.~Biscani, dcgp: Differentiable cartesian genetic programming made easy.
\newblock {\it Journal of Open Source Software\/} {\bf 5}, 2290 (2020).

\bibitem{pontryagin}
L.~Pontryagin, {\it Mathematical Theory of Optimal Processes\/} (CRC Press, 1987).

\bibitem{acciarini2024nonlinear}
G.~Acciarini, N.~Baresi, D.~J. Lloyd, D.~Izzo, Nonlinear propagation of non-gaussian uncertainties.
\newblock {\it Journal of Guidance, Control, and Dynamics\/} {\bf 48}, 903--913 (2024).

\bibitem{michel2022esa}
P.~Michel, M.~K{\"u}ppers, A.~C. Bagatin, B.~Carry, S.~Charnoz, J.~De~Leon, A.~Fitzsimmons, P.~Gordo, S.~F. Green, A.~H{\'e}rique, {\it et~al.\/}, The esa hera mission: detailed characterization of the dart impact outcome and of the binary asteroid (65803) didymos.
\newblock {\it The planetary science journal\/} {\bf 3}, 160 (2022).

\bibitem{lotka1925elements}
A.~J. Lotka, {\it Elements of physical biology\/} (Williams \& Wilkins, 1925).

\bibitem{NODE4LV2021}
W.~Bonnaffé, B.~C. Sheldon, T.~Coulson, Neural ordinary differential equations for ecological and evolutionary time-series analysis.
\newblock {\it Methods in Ecology and Evolution\/} {\bf 12}, 1301-1315 (2021).

\bibitem{richardson2006binary}
D.~C. Richardson, K.~J. Walsh, Binary minor planets.
\newblock {\it Annu. Rev. Earth Planet. Sci.\/} {\bf 34}, 47--81 (2006).

\bibitem{rivkin2021double}
A.~S. Rivkin, N.~L. Chabot, A.~M. Stickle, C.~A. Thomas, D.~C. Richardson, O.~Barnouin, E.~G. Fahnestock, C.~M. Ernst, A.~F. Cheng, S.~Chesley, {\it et~al.\/}, The double asteroid redirection test (dart): planetary defense investigations and requirements.
\newblock {\it The Planetary Science Journal\/} {\bf 2}, 173 (2021).

\bibitem{poincare1892methodes}
H.~Poincar{\'e}, {\it Les m{\'e}thodes nouvelles de la m{\'e}canique c{\'e}leste\/}, vol.~1 (Gauthier-Villars et fils, Paris, 1892). Solutions p{\'e}riodiques. Non-existence des int{\'e}grales uniformes. Solutions asymptotiques.

\bibitem{hamiltonianneuralodes}
S.~Greydanus, M.~Dzamba, J.~Yosinski, {\it NeurIPS\/}, H.~M. Wallach, H.~Larochelle, A.~Beygelzimer, F.~d'Alché Buc, E.~B. Fox, R.~Garnett, eds. (2019), pp. 15353--15363.

\bibitem{goldberg2019juventas}
H.~R. Goldberg, {\"O}.~Karatekin, B.~Ritter, A.~Herique, P.~Tortora, C.~Prioroc, B.~G. Gutierrez, P.~Martino, I.~Carnelli, {\it Small Satellite Conference\/} (2019).

\bibitem{scheeres2013design}
D.~Scheeres, B.~Sutter, A.~Rosengren, Design, dynamics and stability of the osiris-rex sun-terminator orbits.
\newblock {\it Advances in the Astronautical Sciences\/} {\bf 148}, 3263--3282 (2013).

\bibitem{raducan2024physical}
S.~Raducan, M.~Jutzi, A.~Cheng, Y.~Zhang, O.~Barnouin, G.~Collins, R.~Daly, T.~Davison, C.~Ernst, T.~Farnham, {\it et~al.\/}, Physical properties of asteroid dimorphos as derived from the dart impact.
\newblock {\it Nature Astronomy\/} {\bf 8}, 445--455 (2024).

\bibitem{glassmeier2007rosetta}
K.-H. Glassmeier, H.~Boehnhardt, D.~Koschny, E.~K{\"u}hrt, I.~Richter, The rosetta mission: flying towards the origin of the solar system.
\newblock {\it Space Science Reviews\/} {\bf 128}, 1--21 (2007).

\bibitem{DNN_improved_tracking}
Q.~Li, J.~Qian, Z.~Zhu, X.~Bao, M.~K. Helwa, A.~P. Schoellig, {\it 2017 IEEE International Conference on Robotics and Automation (ICRA)\/} (2017), pp. 5183--5189.

\bibitem{osti_10100589}
G.~Tang, W.~Sun, K.~Hauser, Learning trajectories for real- time optimal control of quadrotors.
\newblock {\it IEEE/RSJ Intl Conf on Intelligent Robots and Systems\/}  (2018).

\bibitem{Kaufmann2023}
E.~Kaufmann, L.~Bauersfeld, A.~Loquercio, M.~M{\"u}ller, V.~Koltun, D.~Scaramuzza, Champion-level drone racing using deep reinforcement learning.
\newblock {\it Nature\/} {\bf 620}, 982-987 (2023).

\bibitem{aggressive_online_control}
S.~Li, E.~{\"{O}}zt{\"{u}}rk, C.~D. Wagter, G.~C. H.~E. de~Croon, D.~Izzo, Aggressive online control of a quadrotor via deep network representations of optimality principles.
\newblock {\it CoRR\/} {\bf abs/1912.07067} (2019).

\bibitem{cheng2018real}
L.~Cheng, Z.~Wang, F.~Jiang, C.~Zhou, Real-time optimal control for spacecraft orbit transfer via multiscale deep neural networks.
\newblock {\it IEEE Transactions on Aerospace and Electronic Systems\/} {\bf 55}, 2436--2450 (2018).

\bibitem{federici2021deep}
L.~Federici, B.~Benedikter, A.~Zavoli, Deep learning techniques for autonomous spacecraft guidance during proximity operations.
\newblock {\it Journal of Spacecraft and Rockets\/} {\bf 58}, 1774-1785 (2021).

\bibitem{dario_seb_gcnet}
D.~Izzo, S.~Origer, Neural representation of a time optimal, constant acceleration rendezvous.
\newblock {\it Acta Astronautica\/} {\bf 204}, 510-517 (2023).

\bibitem{origer2023quad}
S.~Origer, C.~D. Wagter, R.~Ferede, G.~C. H.~E. de~Croon, D.~Izzo, Guidance \& control networks for time-optimal quadcopter flight (2023).

\bibitem{origer2024neuralODEfix}
S.~Origer, D.~Izzo, Closing the gap: Optimizing guidance and control networks through neural odes (2024).

\bibitem{izzo2021real}
D.~Izzo, E.~{\"O}zt{\"u}rk, Real-time guidance for low-thrust transfers using deep neural networks.
\newblock {\it Journal of Guidance, Control, and Dynamics\/} {\bf 44}, 315--327 (2021).

\bibitem{sitzmann2020implicit}
V.~Sitzmann, J.~Martel, A.~Bergman, D.~Lindell, G.~Wetzstein, Implicit neural representations with periodic activation functions.
\newblock {\it Advances in neural information processing systems\/} {\bf 33}, 7462--7473 (2020).

\bibitem{origer2024gcnetperiodic}
S.~Origer, D.~Izzo, Guidance and control networks with periodic activation functions (2024).

\bibitem{wittig2015propagation}
A.~Wittig, P.~Di~Lizia, R.~Armellin, K.~Makino, F.~Bernelli-Zazzera, M.~Berz, Propagation of large uncertainty sets in orbital dynamics by automatic domain splitting.
\newblock {\it Celestial Mechanics and Dynamical Astronomy\/} {\bf 122}, 239--261 (2015).

\bibitem{berz1998verified}
M.~Berz, K.~Makino, Verified integration of odes and flows using differential algebraic methods on high-order taylor models.
\newblock {\it Reliable computing\/} {\bf 4}, 361--369 (1998).

\bibitem{biscani2022reliable}
F.~Biscani, D.~Izzo, Reliable event detection for taylor methods in astrodynamics.
\newblock {\it Monthly Notices of the Royal Astronomical Society\/} {\bf 513}, 4833--4844 (2022).

\bibitem{biscani2025heyoka}
F.~Biscani, D.~Izzo, {heyoka: High-precision Taylor integration of ordinary differential equations} (2025). Accessed: 2025-03-23.

\bibitem{valli2013nonlinear}
M.~Valli, R.~Armellin, P.~Di~Lizia, M.~R. Lavagna, Nonlinear mapping of uncertainties in celestial mechanics.
\newblock {\it Journal of Guidance, Control, and Dynamics\/} {\bf 36}, 48--63 (2013).

\bibitem{lattner2004llvm}
C.~Lattner, V.~Adve, {\it International symposium on code generation and optimization, 2004. CGO 2004.\/} (IEEE, 2004), pp. 75--86.

\bibitem{krantz2002implicit}
S.~G. Krantz, H.~R. Parks, {\it The implicit function theorem: history, theory, and applications\/} (Springer Science \& Business Media, 2002).

\bibitem{izzo2022geodesy}
D.~Izzo, P.~G{\'o}mez, Geodesy of irregular small bodies via neural density fields.
\newblock {\it Communications Engineering\/} {\bf 1}, 48 (2022).

\bibitem{torchdyn}
M.~Poli, S.~Massaroli, A.~Yamashita, H.~Asama, J.~Park, Torchdyn: A neural differential equations library.
\newblock {\it arXiv preprint arXiv:2009.09346\/}  (2020).

\bibitem{torchdiffeq}
R.~T.~Q. Chen, torchdiffeq (2018).

\end{thebibliography}

\bibliographystyle{ScienceAdvances}

\section*{Acknowledgments}
The authors express their gratitude to Alexander Hadjiivanov, Guido de Croon, and Sean Cowan for their valuable feedback and insightful discussions on the manuscript prior to its submission.

\section*{Authors Contribution}
D.I., S.O., and G.A. conceived the study, developed the theoretical framework, and designed the experiments. D.I. coded some of the preliminary tools and wrote the initial draft of the manuscript, while S.O. and G.A. coded and performed the experiments. F.B. developed the code base for creating and manipulating any-order variational equations and their efficient numerical integration. All authors reviewed the manuscript, analyzed the results, and contributed to the final version.

\section*{Competing interests} 
The authors declare that they have no competing interests. 

\end{document}